\documentclass[12pt]{amsart}
\begin{document}



\def\ron{\noindent{\bf Ron:\ }}
\def\endron{ \hfill\rule{10mm}{.75mm} \break}

\def\tom{\noindent{\bf Tom:\ }}
\def\endtom{ \hfill\rule{10mm}{.75mm} \break}

\def\jim{\noindent{\bf Jim:\ }}
\def\endjim{ \hfill\rule{10mm}{.75mm} \break}

\def\VV{\mathbb{L}}
\def\gg{$\mathfrak g$}
\def\calA{\mathcal A}
\def\cal{\mathcal}
\def\be{\begin{eqnarray}}

\def\ee{\end{eqnarray}}
\def\Amod{$\calA$-module}
\def\G{\Gamma}

\def\D{\Delta}

\def\qed{\hbox{${\vcenter{\vbox{

\hrule height 0.4pt\hbox{\vrule width 0.4pt height 6pt

\kern5pt\vrule width 0.4pt}\hrule height 0.4pt}}}$}}

\newtheorem{theorem}{Theorem}

\newtheorem{thm}{Theorem}

\newtheorem{lemma}{Lemma}

\newtheorem{definition}{Definition}

\newtheorem{corollary}{Corollary}
\newtheorem{prop}{Proposition}

\def\prl{\partial}
\def\XX{\Xi}
\def\PP{\Phi}
\def\LP{\Lambda^*\Phi}
\def\DD{\Delta}

\def\CD{{\rm Coder}(\Lambda^*\Phi)}
\def\TT{\Theta}
\def\LPP{{\rm Hom}(\LP,\PP)}
\def\LPX{{\rm Hom}(\LP,\XX)}
\def\dd{\delta}
\def\ddd{\hat\delta}
\def\ddf{\widebar{\dd(f)}}
\def\ddh{\widebar{\dd(h)}}
\def\ddg{\widebar{\dd(g)}}
\def\widebar{\overline}
\def\longoverbar{\overline}
\def\ocirc{\odot}
\def\hdd{\hat\delta}
\def\LV{\Lambda^* (\downarrow \VV)}
\def\LVV{\Lambda^* (\downarrow \VV)}

\def\ww{\wedge}
\def\BBvD{Berends, Burgers and van Dam }
\def\BBvDc{Berends, Burgers and van Dam,}
\def\pp{\phi} \def\ppi{\pp_1\wedge \cdots \wedge \pp_ } \def\ss{\sigma}

\pagestyle{myheadings}
\markboth{Fulp, Lada, Stasheff}{Sh-Lie Gauge Transformations}

\title[Sh-Lie algebras Induced by Gauge Transformations] {Sh-Lie algebras
Induced by Gauge Transformations}

\author{Ron Fulp}
\address{Department of Mathematics, North Carolina State University,
Raleigh NC 27695}
\email{fulp@math.ncsu.edu}
\author{Tom Lada}
\address{Department of Mathematics, North Carolina State University,
Raleigh NC 27695}
\email{lada@math.ncsu.edu}
\author{Jim Stasheff }

\address{Department of Mathematics, University of North Carolina, Chapel
Hill, NC 27599-3250, USA}
\email{jds@math.unc.edu}
\thanks{Stasheff's research supported in part by the NSF throughout most of his
career, most recently under grant DMS-9803435. }





\def\gg{$\mathfrak g$}
\def\cal{\mathcal}
\def\be{\begin{eqnarray}}

\def\ee{\end{eqnarray}}

\def\G{\Gamma}

\def\D{\Delta}

\def\qed{\hbox{${\vcenter{\vbox{

\hrule height 0.4pt\hbox{\vrule width 0.4pt height 6pt

\kern5pt\vrule width 0.4pt}\hrule height 0.4pt}}}$}}

\def\prl{\partial}
\def\XX{\Xi}
\def\PP{\Phi}
\def\LP{\Lambda^*\Phi}
\def\DD{\Delta}

\def\CD{{\rm Coder}(\Lambda^*\Phi)}
\def\TT{\Theta}
\def\LPP{{\rm Hom}(\LP,\PP)}
\def\LPX{{\rm Hom}_k(\LP,\XX)}
\def\dd{\delta}
\def\ddd{\hat\delta}
\def\ddf{\widebar{\dd(f)}}
\def\ddh{\widebar{\dd(h)}}
\def\ddg{\widebar{\dd(g)}}
\def\widebar{\overline}
\def\longoverbar{\overline}
\def\ocirc{\odot}
\def\hdd{\hat\delta}
\def\LV{\Lambda^* (\downarrow V)}
\def\ww{\wedge}
\def\ppi{\pp_1\wedge \cdots \wedge \pp_ } \def\ss{\sigma}

\subjclass{}

\begin{abstract}
Traditionally symmetries of
field theories are encoded via Lie group actions, or more generally, as Lie
algebra actions. A significant generalization is required when `gauge
parameters' act in a field dependent way. Such symmetries
appear in several field theories, most
notably in a `Poisson induced' class due to Schaller and Strobl
\cite{schaller-strobl} and to Ikeda\cite{ikeda}, and employed by Cattaneo
and Felder
\cite{cf:defquant} to implement Kontsevich's deformation quantization
\cite{kont:defquant}.
Consideration of `particles of spin $> 2$ led \BBvD
\cite{burgers:diss,BBvd:three,BBvD:probs} to study `field dependent
parameters' in a
setting permitting an analysis in terms of smooth functions.
Having recognized the resulting structure as that of an sh-lie algebra
($L_\infty$-algebra), we have now formulated such structures entirely
algebraically
and applied it to a more general class of theories with field dependent
symmetries.
\end{abstract}

\maketitle

\section{Introduction}

Ever since the discovery of Yang-Mills theory, physicists have been
intrigued by the different
manifestations of symmetries in field theories. Symmetries in gravitational
theories are induced by
spacetime transformations which preserve the spacetime structure whereas
Yang-Mills symmetries are defined
via transformations of some internal vector space. Many authors have
attempted to reformulate
gravitational symmetries in a manner which is compatible with the
Yang-Mills approach as quantization of
Yang-Mills theories is better understood than most attempts to quantize
gravity.

The present paper has as its purpose to show that gauge symmetries of
certain field theories have an
unexpectedly rich algebraic structure. Traditional theories lead one to
expect that the symmetries of
field theories are encoded via Lie group actions, or more generally, as Lie
algebra actions. We find that
the gauge symmetries of many field theories in fact do not arise from a Lie
algebra action, but rather from an sh-Lie (or $L_\infty$) algebra action.

The physics of ``particles of spin $\leq 2$'' leads to representations of a
Lie algebra $\Xi$ of
gauge parameters on a vector space $\Phi$ of fields. A significant
generalization occurs when the gauge parameters act in a {\em field
dependent} way.
By a field dependent action of $\XX$ on $\PP,$ \BBvD
\cite{burgers:diss,BBvD:arb,BBvD:probs} mean a polynomial (or power series)
map $\dd (\xi)(\pp) =
\Sigma_{i\geq 0} T_i(\xi,\pp)$ where
$T_i$ is linear in $\xi$ and polynomial of homogeneous degree $i$ in $\pp.$

Field dependent gauge symmetries appear in several field theories, most
notably in a `Poisson induced' class due to Schaller and Strobl
\cite{schaller-strobl}
and to Ikeda \cite{ikeda}, and
employed by Cattaneo and Felder
\cite{cf:defquant} to implement Kontsevich's deformation quantization
\cite{kont:defquant}. Ikeda
\cite{ikeda} considers two-dimensional and three-dimensional \cite{ikeda:3d}
theories with a generalized
Yang-Mills field which has values in a
so-called nonlinear Lie algebra. He finds that if the non-linear Lie
structure is chosen appropriately and
if he allows the Yang-Mills field to interact with certain scalar fields,
then he can recapture gravitational
theories in two dimensions. In this way, two-dimensional gravity is
formulated as a Yang-Mills theory and
its symmetries arise in the same way as traditional Yang-Mills symmetries.
The three-dimensional case \cite{ikeda:3d} provides deformations of
physicists' BF theories and analogous results hold in higher
dimensions.

Although expressed rather differently, the \BBvD approach provides further
insight into
the algebraic structure of the gauge symmetries of the above class of field
theories. In fact their context is more general than that of Ikeda and that
of Cattaneo and Felder, since \BBvD consider arbitrary field theories,
subject only to the requirement that the commutator of two gauge symmetries
be another gauge symmetry whose gauge parameter is possibly field dependent.
We refer to this requirement as the BBvD hypothesis. Notice
\BBvD do not require an a
priori given
Lie structure to induce the algebraic structure of the gauge symmetry
``algebra". On the other hand,
Ikeda requires a
structure called a nonlinear Lie algebra which he uses to obtain symmetries
which in turn are used to find
a Lagrangian for which the symmetries are gauge symmetries. In this sense,
his nonlinear Lie structure
drives the entire theory. Similarly, Cattaneo and Felder have a Poisson
structure which explicitly appears
in both the action of their theory and in their gauge symmetries.

The present work has as its goal to clarify the algebraic structure of the
more general gauge ``algebra" outlined in \BBvD.
When the BBvD hypothesis is satisfied, we show that the gauge symmetry
algebra of a
large class of field theories is an sh-Lie algebra. Of course, as we show,
this sh-Lie structure, in special cases, will reduce to the more familiar
Lie structures one encounters in various field theories. On the other hand,
some of these field theories satisfy the BBvD hypothesis
only `on-shell'. When closure on the
original space of parameters is lost, physicists speak of an `open algebra'.
This leads us, in section 7, to a `generalized BBvD hypothesis'
which in turn will allow us to show how the sh-Lie structure must be
modified to handle `off-shell' gauge symmetries.

We formulate the relevant structures in BBvD's theory in terms of linear
maps from a
certain coalgebra $\Lambda^*\Phi$ into the respective vector spaces $\Phi$ of
fields and $\Xi$ of
gauge parameters. The coalgebra and the algebra structures of $\Lambda^*\Phi$
as well as the Lie algebra structure of $\LPP$ are described  in Section 2.
It turns out that the space $\Xi$ of gauge
parameters has, in general, no natural Lie structure, but the space of
linear maps from
$\Lambda^*\Phi$ into $\Xi$ is a Lie algebra under certain mild
assumptions along with 
the BBvD hypothesis. This is proved in Section 3. Section 4 provides the
reader with a
short description of two equivalent methods for defining sh-Lie algebras.
Our main
result is found in Section 5 where we show that, under the same assumptions
required
in Section 3, the fields and gauge parameters combine to form an sh-Lie
algebra.
In Section 6 we show how our results
relate to the
classical situation in which the space $\Xi$ of gauge parameters is a Lie
algebra which
acts on the space $\Phi$ of fields. Section 7 provides further links to the
physics
literature where certain sigma-models are known to satisfy the BBvD
hypothesis only `on-shell'. This requires us to further generalize the BBvD
hypothesis;
consequently these gauge algebras are ``on shell" sh-Lie algebras which are not
``on shell" Lie algebras. Finally, in Section 8 we show explicitly how our
formalism
applies to the work of
Ikeda \cite{ikeda} on two-dimensional gravitational theories and his study
of non-linear Lie algebras.
In addition, we show that Ikeda's bracket is the `non-linear' analog of the
Kirillov-Kostant bracket.

We are grateful to \BBvD ~for the inspiration of Burgers' dissertation and
especially to van Dam for several discussions as our research developed.

\section{Our framework}

We work with vector spaces over a field $k$ of characteristic 0 or, more
generally, over
a commutative $k$-algebra $\calA,$ typically, $C^\infty (M)$ for some
smooth manifold $M.$
Unless otherwise specified, $Hom$ will denote the $\calA$-module of
$\calA$-linear maps.

Let $\PP$ be a free $\calA$-module and let $\LP$ denote the
free nilpotent graded cocommutative coalgebra over $\calA$ cogenerated by
$\PP$ with comultiplication
denoted $\DD.$ This is the coalgebra
of graded symmetric tensors in the full tensor coalgebra on $\PP$. The
$\calA$-module $\CD$
of {\em coderivations} (over $\calA$) on $\LP$ is a Lie algebra with
bracket given by the commutator with
respect to composition. Recall that a coderivation is a linear map
$\TT:\LP\to\LP$ that satisfies the equation
$$ \DD\circ\TT = (\TT\otimes 1 + 1\otimes \TT)
\circ \DD.$$ (In the graded situation, the usual Koszul sign conventions
are in effect.)

The $\calA$-module $\LPP$ is isomorphic to $\CD$ and hence inherits a Lie
algebra structure; the bracket
on $\LPP$ is known as the Gerstenhaber bracket
\cite{gerst:coh,jds:intrinsic}. The isomorphism $$ \LPP\ni h
\rightleftharpoons \bar h \in \CD $$ is given by the correspondence $$
\bar h (\pp_1\wedge \cdots \wedge \pp_n) = \underset{\{unshuff\}}{\Sigma}
h(\pp_{\ss(1)}\wedge \cdots
\wedge \pp_{\ss(p)}) \wedge
\pp_{\ss(p+1)}\wedge\cdots\wedge \pp_{\ss(n)}$$ for $h\in {\rm
Hom}(\Lambda^p(\PP),\PP).$ The set
$\{unshuff\}$ is the set of $(p,n-p)$-{\em
unshuffles},
that is, the permutations of $\{1,\dots,n\}$ such that $\sigma (1)<\cdots
<\sigma(p)$ and
$\sigma(p+1)<\cdots<\sigma(n)$.
We may write
$\bar h$ as the composition
$\bar h=m\circ(h\otimes1)\circ\Delta$ where $m$ is the usual product in
$\Lambda^*\Phi$
regarded as an algebra (symmetric on even elements and skew on odd ones; no
compatability with the
coproduct is assumed nor needed).

The Gerstenhaber bracket on $Hom(\Lambda^*\Phi,\Phi)$ may be described as
$[f,g] =
f\circ\bar g - g\circ \bar f$
where $\bar f$ and $\bar g$ are the coderivations corresponding to $f$ and
$g.$ In this notation
the ``Gerstenhaber comp" operation may be defined by $f\odot g= f\circ\bar
g,$ for $f,g\in
Hom(\Lambda^*\Phi,\Phi).$
Thus an alternative notation for the Lie bracket on
$Hom(\Lambda^*\Phi,\Phi)$ is
$[f,g] = f\odot g - g\odot f.$

\section{A preliminary result}

Now let $\XX$ and $\PP$ be arbitrary \Amod s. In the Yang Mills example,
the map
$\delta$ takes gauge
parameters to covariant derivatives. In generalizing that, we suppose that
we are
given a $k$-linear map $\dd:\XX\to \LPP.$ Formally, we can
write $\dd(\xi) = \Sigma_{i=0} T_i(\xi)$ where $T_i$ is 0 except on
$\Lambda^i\PP$.
(This $T_i$ is equivalent to the $T_i$ of \BBvD.)
We extend $\delta$ to a map $$\hat\delta: \LPX \to \LPP$$ by
$$\hat\dd (\pi) = ev\circ(\dd\circ\pi \otimes 1)\circ \DD$$ where $ev$ is
the evaluation map. That is,

$$\hat\dd (\pi)(\ppi n) = \underset{\{unshuff\}}{\Sigma} \delta
(\pi(\pp_{\ss(1)}\wedge \cdots \wedge \pp_{\ss(p)}) (\pp_{\ss(p+1)}
\wedge\cdots\wedge \pp_{\ss(n)}).$$

We may think of $\XX$ as being contained in $\LPX$ by identifying
$\xi\in\XX$ with the map, also denoted $\xi$, in $\LPX$ which is 0 except
on the scalars
where $\xi(1) =\xi.$ Note that $\wedge^*\Phi$ is an ${\cal A}$ module and
$k\subset {\cal A}$ and so
$1\in k\subset {\cal A}.$ We will be careful to distinguish $k$-linear maps
from ${\cal A}$-linear
as the need occurs. It is easy to see that $\hat\dd(\xi) = \dd(\xi).$

Our problem concerns possible algebraic structure on $\XX;$ consequently we
consider the possibility of
constructing a Lie-type bracket on
\newline $\LPX$ via the mapping $\hat\delta.$ Under certain conditions,
 such a bracket may then be used
to obtain a bracket on the parameter space defined by restricting the
induced bracket on
$\LPX$ to the parameter space $\XX$. With this in mind, define
$$[\pi_1,\pi_2]:=\pi_1\circ
\overline {\hat\delta (\pi_2)}-\pi_2\circ \overline {\hat\delta (\pi_1)},$$ for
$\pi_1,\pi_2\in \LPX.$ It turns out that this bracket does not generally
satisfy the Jacobi
identity. Moreover, if we choose
$\pi_1=\xi,\pi_2=\eta
\in
\XX,$ then
$$
[\xi,\eta] = \xi\odot \hdd (\eta) -\eta\odot\hdd (\xi) = 0, $$ and as a
result, the restriction of the
induced bracket to
$\XX$ yields an abelian Lie algebra
structure. In many cases of interest, the parameter space has an a priori
nonabelian Lie algebra
structure on it and we would
certainly want the Lie structure
on $\LPX$
to reproduce this structure when restricted to the parameter space $\XX.$

In order to assure the Jacobi property of bracket on $\LPX,$ we introduce a
{\em correction term.} We accomplish this, following Berends, Burgers and
van Dam, by assuming that there is a map $$ C: \XX\otimes \XX \to \LPX $$
such that
$$[\dd (\xi), \dd (\eta)] = \hat \dd C(\xi,\eta) \in \LPP$$ for all
$\xi,\eta\in
\XX.$
We will refer to this
as the {\em BBvD hypothesis}.

Extend $ C$ to a mapping
$$\hat C:\LPX\otimes \LPX\to\LPX$$
by
$$
\hat C(\pi_1,\pi_2) = C\circ ((\pi_1\otimes \pi_2)\otimes 1)\circ
(\Delta\otimes 1)\circ \Delta,$$
where we have identified $C$ with its adjoint mapping, which is the mapping
from $\XX \otimes \XX
\otimes \LP$ into $\XX$ defined by $$(\xi,\eta,\phi_1\wedge\cdots\wedge
\phi_n)\longrightarrow
C(\xi,\eta)(\phi_1\wedge\cdots\wedge \phi_n).$$

Next, we redefine the bracket on $\LPX$ given above by including the
correction term $C:$
$$
[\pi_1,\pi_2] := \pi_1\odot \hat\dd(\pi_2) - \pi_2\odot \hat\dd(\pi_1) +
\hat C(\pi_1,\pi_2).
$$
\begin{thm} The mapping $\hat\delta$ preserves brackets; that is, $\hat\dd
[\pi_1,\pi_2] =[\hat\dd(\pi_1),\hat\dd(\pi_2)].$ Moreover, if
$\hat\dd:\LPX\to \LPP$ is injective, then $[\pi_1,\pi_2]$ satisfies the
Jacobi identity. \end{thm}

\begin{proof}
Observe that if $\pi_1,\pi_2\in \LPX,$
\begin{align*}\hat\delta(\pi_1)\odot
\hat\delta(\pi_2)=
& \hat\delta(\pi_1)\circ
\longoverbar{\hat\delta(\pi_2)}\\
& = ev\circ [(\hat\delta\circ
\pi_1)\otimes
1]\circ \Delta \circ \longoverbar{\hat\delta(\pi_2)} \\
&=ev\circ [(\hat\delta\circ \pi_1)\otimes 1]\circ
\{\longoverbar{(\hat\delta(\pi_2)}\otimes 1) +(1\otimes
\longoverbar{\hat\delta(\pi_2)})\circ \Delta\}\\ & =ev \circ
[((\hat\delta\circ \pi_1)\circ
\longoverbar {\hat\delta(\pi_2)})\otimes 1]\circ \Delta+ev \circ
[((\hat\delta\circ \pi_1)\otimes
\longoverbar {\hat\delta(\pi_2)})]\circ \Delta\\ & =\hat\delta
(\pi_1\circ\longoverbar{\hat\delta(\pi_2)}))+ ev \circ [((\hat\delta\circ
\pi_1)\otimes \longoverbar {\hat\delta(\pi_2)})]\circ \Delta\\& =\hat\delta
(\pi_1\odot \hat\delta(\pi_2)))+ ev \circ [((\hat\delta\circ \pi_1)\otimes
\longoverbar {\hat\delta(\pi_2)})]\circ \Delta.\end{align*}

It follows that
$$[\hat\delta(\pi_1),\hat\delta(\pi_2)]=
\hat\delta (\pi_1\odot \hat\delta(\pi_2)))-\hat\delta (\pi_2\odot
\hat\delta(\pi_1)))+E$$ where
$$E=ev \circ
[((\hat\delta\circ \pi_1)\otimes \longoverbar { \hat\delta(\pi_2)})]\circ
\Delta -ev \circ [(\hat\delta\circ
\pi_2)\otimes \longoverbar { \hat\delta(\pi_1)}]\circ \Delta.$$ This says that
$E$ measures the deviation of $\hat\delta$ from being a $\CD$-module map.

We must show that E is in the image of $\hat\delta.$ Recall that for $f\in
\LPP,$ we have
$\longoverbar {f}=m\circ (f\otimes 1)\circ \Delta,$ where $m$ denotes the
algebra
(wedge) product on $\LP.$ Thus
$$(\hat\delta\circ \pi_1)\otimes
\longoverbar{\hat\delta(\pi_2)}=(\hat\delta\circ \pi_1) \otimes \{m\circ
(\hat\delta(\pi_2)\otimes 1)\circ \Delta\}$$ $$=(\hat\delta
\circ\pi_1)\otimes \{m\circ ([ev\circ ((\hat\delta\circ \pi_2)\otimes
1)]\otimes 1)\circ (\Delta\otimes 1)\circ \Delta\}.$$ $$=(\hat\delta
\circ\pi_1)\otimes \{m\circ ([ev\circ ((\hat\delta\circ \pi_2)\otimes
1)]\otimes 1)\circ (1\otimes \Delta)\circ \Delta\}.$$

For $F\in \LP,$ write $\Delta(F)=\sum(F_1\otimes F_2), \Delta
(F_2)=\sum(F_{21}\otimes F_{22})$
and $\Delta(F_{22})=\sum(F_{221}\otimes F_{222}).$ In order to simplify
notation we drop the
summation symbol wherever the latter coproducts appear below. From our last
calculation we have
$$(ev \circ [(\hat\delta\circ \pi_1)\otimes \longoverbar {
\hat\delta(\pi_2)}]\circ \Delta)(F)= [\hat\delta(\pi_1(F_1))\odot
\hat\delta(\pi_2(F_{21}))](F_{22}),$$ and $$(ev \circ [(\hat\delta\circ
\pi_2)\otimes \longoverbar {\Delta(\pi_1)}]\circ \Delta)(F)=
[\hat\delta(\pi_2(F_1))\odot \hat\delta(\pi_1(F_{21}))](F_{22}).$$ Because
$\Delta$ is cocommutative, the full summations are equal:
$$\Sigma F_1\otimes
F_{21} \otimes F_{22} = \Sigma F_{21}\otimes F_1 \otimes F_{22}.$$ Thus
$$E(F)=[\hat\delta(\pi_1(F_1)),
\hat\delta(\pi_2(F_{21}))](F_{22})=\hat\delta(C(\pi_1(F_1),\pi_2(F_{21}),F_{221}
))(F _{2 22})$$
$$=ev \circ(\{\hat\delta \circ C \circ [(\pi_1\otimes \pi_2)\otimes
1]\}\otimes 1)(F_1\otimes F_{21}\otimes
F_{221}\otimes F_{222})$$
$$=ev \circ (\{\hat\delta \circ C \circ [(\pi_1\otimes \pi_2)\otimes
1]\}\otimes 1)(F_1\otimes F_{21}\otimes
(\Delta F_{22}))$$
$$=ev \circ (\{\hat\delta \circ C \circ [(\pi_1\otimes \pi_2)\otimes
1]\}\otimes 1)
( ([1\otimes((1\otimes\Delta)\circ \Delta)]\circ \Delta)(F))$$ $$=ev \circ
(\{\hat\delta \circ C \circ
[(\pi_1\otimes \pi_2)\otimes 1]\}\otimes 1) ((1\otimes 1\otimes
\Delta)\circ(1\otimes \Delta)\circ \Delta)(F)).$$ It follows from
coassociativity that $$E=ev \circ (\{\hat\delta \circ C \circ
[(\pi_1\otimes \pi_2)\otimes 1]\}\otimes 1)\circ ((\Delta\otimes 1\otimes
1)\circ(\Delta\otimes 1)\circ\Delta $$ 
$$=ev \circ (\{\hat\delta \circ
C \circ [(\pi_1\otimes \pi_2)\otimes 1]\}\otimes 1)\circ ([(\Delta \otimes
1)\circ \Delta]\otimes 1)\circ \Delta $$ $$=ev \circ (\{\hat\delta \circ C
\circ
[(\pi_1\otimes \pi_2)\otimes 1]\circ (\Delta \otimes 1)\circ \Delta
\}\otimes 1)\circ \Delta $$
$$=ev \circ (\{\hat\delta \circ \hat C(\pi_1,\pi_2)\}\otimes 1)\circ
\Delta= \hat\delta(\hat C(\pi_1,\pi_2)).$$ Thus $E$ is in the image of
$\hat\delta$ and in fact $$[\hat\delta(\pi_1),\hat\delta
(\pi_2)]=\hat\delta (\pi_1\odot \hat\delta(\pi_2)-\pi_2\odot
\hat\delta(\pi_1)) +\hat\delta (\hat
C(\pi_1,\pi_2))=\hat\delta([\pi_1,\pi_2]).$$

To verify the Jacobi identity, apply $\hat\dd$ to the Jacobi expression in
$\LPX.$ By the morphism
condition just established, the result is the Jacobi identity valid in
$\LPP.$ Assuming that $\hat\dd$ is injective, the Jacobi identity in $\LPX$
follows.
\end{proof}

This result suggests that the parameter space should be enlarged to include
all of
$Hom_k(\Lambda^*\Phi,\Xi).$ It
turns out that the polynomial equations of physical relevance define an
sh-Lie structure
on an appropriate graded vector space $\VV.$ We consider the sh-Lie formalism
briefly in the next section.

\section{Sh-Lie algebras}
We now review the relationship between sh-Lie algebras
($L_\infty$-algebras) and cocommutative
coalgebras \cite{ls,lada-markl}.
Let $(\VV,d)$ be a differential graded vector space. If $(\VV,d)$ is a chain
complex (degree $d=-1$),
then an sh-Lie structure on $\VV$ is a collection of skew symmetric linear
maps $l_n:\VV^{\otimes n}\longrightarrow \VV$ of degree $n-2$ that satisfy
the relations $$\sum_{i+j=n+1}\sum_{\sigma}
e(\sigma)(-1)^{\sigma}(-1)^{i(j-1)}
l_j(l_i(x_{\sigma(1)},\dots,
x_{\sigma(i)}),\dots ,x_{\sigma(n)}) =0$$ where $(-1)^{\sigma}$ is the sign
of the permutation $\sigma$, $e(\sigma)$ is the sign that arises from the
degrees of the permuted elements and $\sigma$ is taken over all $(i,n-i)$
unshuffles.

If $l_n = 0$ for $n \geq 3$, this is just the description of a dg Lie algebra.

If $(\VV,d)$ is a cochain complex (degree $d=+1$), then the sh-Lie structure
on $\VV$ is given by skew symmetric
linear maps $l_n:\VV^{\otimes n}\longrightarrow \VV$ of degree $2-n$ that
satisfy the same relations.

Let $\uparrow \VV$ denote the suspension of the graded vector space $\VV$; i.e.
$\uparrow \VV$ is the graded
vector space with $(\uparrow \VV)_n=\VV_{n-1}$; similarly, let $\downarrow \VV$
denote the desuspension of $\VV$; i.e.
$(\downarrow \VV)_n=\VV_{n+1}$.

One may then describe an sh-Lie structure on the chain complex $(\VV,d)$ by a
coderivation $\overline D$ of degree
$-1$ on the coalgebra $\Lambda^*(\uparrow \VV)$ such that ${\overline
D}^2=0$; similarly, an sh-Lie structure on the cochain complex $(\VV,d)$ is
a
coderivation $\overline D$ of degree
$+1$ on the
coalgebra $\Lambda^*(\downarrow \VV)$ such that ${\overline D}^2=0$.
Equivalently, the sh-Lie structure may
be described by a linear mapping $D: \Lambda^*(\downarrow \VV)\longrightarrow
(\downarrow \VV)$ such that $D\circ \overline D=0.$ The proof of the
assertion for chain complexes may be found in \cite{ ls} and
\cite{jds:intrinsic}; a proof for cochain complexes can be formulated by a
straightforward modification of the proof for chain complexes.

\section{The gauge algebra is an sh-Lie algebra}

We now restrict our attention to the constant maps in $\LPX$ and show that
our algebraic structure
on $\LPX$ induces an sh-Lie structure on the graded space $\VV=\{\XX, \Phi\}$.
Throughout
this section,
we assume the BBvD hypothesis and that $\hat\delta$ is injective, so
Theorem 1 holds and consequently the bracket on $\LPX$ defined by $$
[\pi_1,\pi_2] := \pi_1\odot \hat\dd(\pi_2) - \pi_2\odot \hat\dd(\pi_1) +
\hat C(\pi_1,\pi_2)
$$
satisfies the Jacobi identity. By definition, $$ [\delta(\xi),\delta(\eta)]
= \dd(\xi)\odot\dd(\eta) - \dd(\eta)\odot\dd(\xi) $$ while the definition
of $ C$ gives $$
[\dd (\xi), \dd (\eta)] = \ddd C(\xi,\eta) \in \LPP, $$ so our {\em
commutator relation} is
$$
\dd(\xi)\odot\dd(\eta) - \dd(\eta)\odot\dd(\xi) = \ddd( C (\xi,\eta)). $$

The definition of the bracket in $\LPX$ restricted to constant maps takes
on the form $[\xi_1,,\xi_2]= C(\xi_1,\xi_2)$. Consequently, the Jacobi identity
takes on the form
$$
[C(\xi_1,\xi_2),\xi_3] -[C(\xi_1,\xi_3),\xi_2] +[C(\xi_2,\xi_3),\xi_1] =0. $$

Let us examine the first term:
\begin{multline*}
[C(\xi_1,\xi_2),\xi_3] = \\
C(\xi_1,\xi_2)\odot \dd(\xi_3) - \xi_3\odot \ddd C(\xi_1,\xi_2) + \hat
C(C(\xi_1,\xi_2),\xi_3) = \\
C(\xi_1,\xi_2)\odot \dd(\xi_3) +\hat C(C(\xi_1,\xi_2),\xi_3)
\end{multline*} because
$\xi_3\odot \ddd C(\xi_1,\xi_2) = 0$ as $\xi$ is a
constant map (non-zero only on scalars). We now add together the results
from
the remaining two terms and write the {\em Jacobi relation} as

$$C(\xi_1,\xi_2)\odot\delta(\xi_3)-C(\xi_1,\xi_3)\odot\delta(\xi_2)
+C(\xi_2,\xi_3)\odot\delta(\xi_1)$$
$$+\hat C(C(\xi_1,\xi_2),\xi_3)-\hat C(C(\xi_1,\xi_3),\xi_2)+\hat
C(C(\xi_2, \xi_3),\xi_1)=0.$$

For the sh-Lie structure, we first combine the fields and gauge parameters
to form a single differential
graded vector space $\VV.$

\begin{definition} The underlying dg vector space $\VV$ of {\em the sh-Lie
algebra} has $\XX$ in
degree 0, $\PP$ in degree 1 and 0 in all other degrees. The differential
$\prl: \XX\to\PP$ is given by $\prl
(\xi) = \dd(\xi)(1)\in \PP.$
\end{definition}

\begin{thm} The linear map $$D:\LVV \to \downarrow \VV $$
given by
\begin{align*}
D(\xi) = &\prl(\xi)\\ D(\xi\wedge\phi_1\ww\cdots\ww\phi_n) =
&\dd(\xi)(\phi_1\wedge\cdots\wedge\phi_n) {\rm \ for \ } n\geq 1 \\
D(\xi_1\wedge\xi_2\wedge\phi_1\wedge\cdots\wedge\phi_n)=
&C(\xi_1,\xi_2)(\phi_1\wedge\cdots\wedge\phi_n) \end{align*}

\noindent and $D=0$ on elements of $\LVV$ with more than two entries from
$\XX$ or with no entry from $\Xi$ gives $\VV$ the structure of an sh-Lie
algebra.
\end{thm}

\noindent{\bf Remark:} Recall that we have assumed as hypothesis for this
theorem that the bracket on
$\LPX$ satisfies the Jacobi identity. According to Theorem 1 this is true
if $\hat\delta$ is
injective. It is not difficult to prove that $\hat\delta$ is injective
whenever $\delta$
is injective.
If we replace the original parameter space with the new parameter space
$\XX/ ker(\delta)$, one has the sh-Lie structure obtained in the proof below.

\begin{proof} We need only evaluate $D\circ \bar D$ on elements of the form
$(\xi_1\wedge\xi_2\wedge\phi_1\wedge\cdots\wedge\phi_n)$ and
$(\xi_1\wedge\xi_2\wedge\xi_3\wedge\phi_1\wedge\cdots\wedge\phi_n)$.

We begin with
$$D\circ\bar D(\xi_1\wedge\xi_2\wedge\phi_1\wedge\cdots\wedge\phi_n)=$$
$$D\{\sum_{\sigma}\delta(\xi_1)(\phi_{\sigma(1)}\wedge\cdots\wedge\phi_{\sigma(i
)})
\wedge\xi_2\wedge\phi_{\sigma(i+1)}\wedge\cdots\wedge\phi_{\sigma(n)}$$
$$-\sum_{\tau}\delta(\xi_2)(\phi_{\tau(1)}\wedge\cdots\wedge\phi_{\tau(j)})
\wedge\xi_1\wedge\phi_{\tau(j+1)}\wedge\cdots\wedge\phi_{\tau(n)}$$
$$\sum_{\rho}C(\xi_1,\xi_2)(\phi_{\rho(1)}
\wedge\cdots\wedge\phi_{\rho(k)})\wedge
\phi_{\rho(k+1)}\wedge\cdots\wedge\phi_{\rho(n)}\}$$ where $\sigma$, $\tau$
and $\rho$ are the evident unshuffles.

This composition is equal to
$$D\{\sum_{\sigma}\xi_2\wedge\delta(\xi_1)(\phi_{\sigma(1)}\wedge\cdots\wedge
\phi_{\sigma(i)})
\wedge\phi_{\sigma(i+1)}\wedge\cdots\wedge\phi_{\sigma(n)}$$
$$-\sum_{\tau}\xi_1\wedge\delta(\xi_2)(\phi_{\tau(1)}\wedge\cdots\wedge
\phi_{\tau(j)})
\wedge\phi_{\tau(j+1)}\wedge\cdots\wedge\phi_{\tau(n)}$$
$$+\sum_{\rho}C(\xi_1,\xi_2)(\phi_{\rho(1)}\wedge\cdots\wedge\phi_{\rho(k)})
\wedge
\phi_{\rho(k+1)}\wedge\cdots\wedge\phi_{\rho(n)}\}$$
$$=\sum_{\sigma}\delta(\xi_2)(\delta(\xi_1)(\phi_{\sigma(1)}\wedge\cdots\wedge
\phi_{\sigma(i)})
\wedge\phi_{\sigma(i+1)}\wedge\cdots\wedge\phi_{\sigma(n)})$$
$$-\sum_{\tau}\delta(\xi_1)(\delta(\xi_2)(\phi_{\tau(1)}\wedge\cdots\wedge
\phi_{\tau(j)})
\wedge\phi_{\tau(j+1)}\wedge\cdots\wedge\phi_{\tau(n)})$$

$$+\sum_{\rho}\delta(C(\xi_1,\xi_2))(\phi_{\rho(1)}\wedge\cdots\wedge
\phi_{\rho(k)})
(\phi_{\rho(k+1)}\wedge\cdots\wedge\phi_{\rho(n)})$$ which is equal to 0 by
the commutator relation.

For the terms of the form
$(\xi_1\wedge\xi_2\wedge\xi_3\wedge\phi_1\wedge\cdots \wedge\phi_n)$, the
only unshuffles that we need to consider are those that result in terms of
the form
$$(\xi_i\wedge\phi_{\sigma(1)}\wedge\cdots\wedge\phi_{\sigma(p)}
\wedge\xi_j\wedge\xi_k
\wedge\phi_{\sigma (p+1)} \wedge\cdots\wedge\phi_{\sigma(n)})\text{ with
}j<k$$ and
$$(\xi_i\wedge\xi_j\wedge\phi_{\tau(1)}\wedge\cdots\wedge\phi_{\tau(q)}\wedge
\xi_k\wedge\phi_{\tau(q+1)}\wedge\cdots\wedge\phi_{\tau(n)})\text{ with
}i<j.$$ Recall that when $i=2$ in the first term and when $j=3, k=2$ in the
second term, a coefficient of $-1$ must be introduced.

So we have
$$D\circ\bar D(\xi_1\wedge\xi_2\wedge\xi_3\wedge\phi_1\wedge\cdots
\wedge\phi_n)=$$
$$D\{\sum_{\sigma}\delta(\xi_i)(\phi_{\sigma(1)}\wedge\cdots\wedge
\phi_{\sigma(p)}\wedge
\xi_j\wedge\xi_k\wedge\phi_{\sigma(p+1)}\wedge\cdots\wedge\phi_{\sigma(n)}$$ $$+
\sum_{\tau}C(\xi_i,\xi_j)(\phi_{\tau(1)}\wedge\cdots\wedge\phi_{\tau(q)})
\wedge\xi_k\wedge\phi_{\tau(q+1)}\wedge\cdots\wedge\phi_{\tau(n)}\}$$

$$=D\{\sum_{\sigma}\xi_j\wedge\xi_k\wedge\delta(\xi_i)(\phi_{\sigma(1)}\wedge
\cdots\wedge\phi_{\sigma(p)})\wedge
\wedge\phi_{\sigma(p+1)}\wedge\cdots\wedge\phi_{\sigma(n)}$$
$$+\sum_{\tau}C(\xi_i,\xi_j)(\phi_{\tau(1)}\wedge\cdots\wedge\phi_{\tau(q)})
\wedge\xi_k\wedge\phi_{\tau(q+1)}\wedge\cdots\wedge\phi_{\tau(n)}\}$$
$$=\sum_{\sigma}C(\xi_i,\xi_j)(\delta(\xi_k)(\phi_{\sigma(1)}\wedge
\cdots\wedge\phi_{\sigma(p)})\wedge
\phi_{\sigma(p+1)}\wedge\cdots\wedge\phi_{\sigma(n)})$$ $$+
\sum_{\tau}C(C(\xi_i,\xi_j)(\phi_{\tau(1)}\wedge\cdots\wedge\phi_{\tau(q)}),
\xi_k(\phi_{\tau(q+1)}\wedge\cdots\wedge\phi_{\tau(n)})$$ which, after
expanding the $i,j,k$ terms of the
unshuffles along with the signs mentioned above, is seen to equal the
Jacobi relation, and hence is
equal to
0. \end{proof}

\section{The classical strict Lie case}

We examine the classical case in which $\Xi$ is a Lie algebra and $\Phi$ is
a Lie module over $\Xi$. Let
us denote the action of $\Xi$ on $\Phi$ by $\xi\cdot\phi.$ We assume that
we have a linear map
$\partial:\Xi\rightarrow\Phi$ that interacts with the Lie module structure
as follows:
$$\partial
[\xi,\eta]_{\Xi}=\xi\cdot(\partial\eta)+\eta\cdot(\partial\xi)$$ where we
have denoted the Lie
bracket on $\Xi$ by $[\cdot,\cdot]_{\Xi}$. As usual the Lie bracket on
$\VV=\Xi\oplus \Phi$ is given by
\begin{equation} [x,y]_{\VV}=
\left\{ \begin{array}{cl}
[x,y]_{\Xi} & \mbox{for $x,y\in \Xi$} \\ x\cdot y & \mbox{for $x\in
\Xi,y\in\Phi$}\\ 0 & \mbox{for $x,y\in \Phi$.}
\end{array} \right. \end{equation}
Similarly denote
the Lie bracket
on $Hom(\Lambda^*\Phi,\Phi)$ by
$[\cdot,\cdot]_{Hom(\Phi)}$ (see Section 1) and on
$Hom_k(\Lambda^*\Phi,\Xi)$ by $[\cdot,\cdot]_{Hom(\Xi)}$ (see Section 2).

Notice that this case is typical of the gauge structure which arises in
fundamental physical theories
such as Yang-Mills theory and basic gravitational theories. For the
Yang-Mills case, the parameter space
$\XX$ is the set of all smooth functions from the space-time M into the Lie
algebra \gg\ of the structure
group G of the theory (for convenience of exposition, we assume that the
principal bundle
of the theory is trivial). The Lie
bracket on the parameter space is the point-wise bracket of two such
parameters. The fields of Yang-Mills
theory are \gg -valued one-forms on M. Note that \BBvD denote the gauge
transformation action of
$\XX$ on $\Phi$ by $\{A,\Lambda\}$ (for $A\in \Phi,\Lambda\in \XX$) rather
than the notation $\Lambda \cdot A$ used above. In this case, this action
is simply the covariant derivative of $\Lambda$ relative to the connection $A.$

Similarly, when the Einstein-Hilbert action is utilized, the parameter
space is the Lie algebra
of all vector fields $\xi$ on the space-time manifold M. Again, in \BBvDc
the background metric
$\eta$ (Minkowski) is presumed and general metrics are written in the form
$\eta+h$ for
an appropriate symmetric tensor $h.$ Thus the fields of the theory are
symmetric tensors
$h.$ The action of a parameter $\xi$ on a field $h$ is the Lie derivative
of $h$ relative
to the vector field $\xi.$ The function $\delta$ is given by $$(\delta
h)_{\mu\nu}=
\partial_{\mu} \xi_{\nu}+\partial_{\nu}\xi_{\mu}+
[(\partial_{\rho}h_{\mu\nu})\xi^{\rho}-h_{\rho\mu}(\partial^{\rho}\xi_{\nu}) - h
_{\rho\nu} (\partial ^{\rho}\xi_{\mu})].$$ \noindent
Details of these two standard examples may be found in Burgers'
dissertation \cite{burgers:diss}.

Notice that using a bracket notation, $[\xi,\phi]_\VV:=\xi\cdot\phi$, for the
action similar to that in
\BBvDc\ the requirement that the bracket be a chain map with respect to
$\partial$ is simply $\partial[\xi,\eta]_\VV=[\xi,\partial\eta]_\VV+
[\eta,\partial\xi]_\VV$. (We already require that $[\cdot,\cdot]_\VV$ restricts
to $[\cdot,\cdot]_{\Xi}$.)

Let us define the "gauge transformation" $ \delta:\Xi\rightarrow
Hom(\Lambda^*\Phi,\Phi)$ by

\begin{equation} \delta(\xi)(\phi)=
\left\{ \begin{array}{cl}
\partial\xi & \mbox{for $\phi=1$} \\
\xi\cdot\phi & \mbox{for $\phi\in\Lambda^1\Phi=\Phi$}\\ 0 & \mbox{for $\phi =
\phi_1\wedge\cdots\wedge\phi_n\in \Lambda^n\Phi$, $n>1$.} \end{array} \right.
\end{equation}

Extend $\delta$ to $\hat\delta:Hom_k(\Lambda^*\Phi,\Xi)\rightarrow
Hom(\Lambda^*\Phi,\Phi)$ by

\begin{equation} \hat\delta(\pi)(\phi)=
\left\{ \begin{array}{cl}
\delta(\pi(\phi_1))(\phi_2)=\partial\pi(\phi) & \mbox{for $\phi_2=1$} \\
\pi(\phi_1)\cdot\phi_2 &
\mbox{for $\phi_2\in\Lambda^1\Phi=\Phi$}\\ 0 & \mbox{otherwise}. \end{array}
\right. \end{equation}
\noindent
Here $1\in k\subset {\cal A}$  while $\phi$ denotes an arbitrary
element of $\Lambda^*\Phi$ and $\Delta(\phi)=\sum\phi_1\otimes \phi_2$.

The canonical bracket on $Hom_k(\Lambda^*\Phi,\Xi)$ that is induced by
$\hat\delta$ and defined below will not
satisfy the
Jacobi identity in general. This bracket is given by
$$[\pi_1,\pi_2]_{Hom(\Xi)}(\phi)=\pi_1\circ\overline{\hat\delta(\pi_2)}(\phi) -\
pi_2\circ\overline
{\hat\delta(\pi_1)}(\phi).$$
Here,
$\overline{\hat\delta(\pi)}(\phi)=
\hat\delta(\pi)(\phi_1)\wedge\phi_2=\delta(\pi(\phi_{11}))(\phi_{12})\wedge
\phi_2,$ and
\begin{equation}
\delta(\pi(\phi_{11}))(\phi_{12})\wedge\phi_2= \left\{ \begin{array}{cl}
\partial(\pi(\phi_1)\wedge\phi_2 & \mbox{if $\phi_{12}=1$}\\
\pi(\phi_{11})\cdot\phi_{12})\wedge\phi_2 & \mbox{if
$\phi_{12}\in\Lambda^1\Phi$}\\
0 & \mbox{otherwise.}
\end{array}
\right.
\end{equation}

In particular, if $\pi(\phi)=\xi(\phi)$ is defined to be the map with value
$\xi$ when $\phi=1\in k$ and $0$ otherwise, then for
$\phi=\sum(\phi_1\wedge \phi_2),$ \begin{equation}
\overline{\hat\delta(\xi)}(\phi)= \left \{ \begin{array}{cl}
\partial\xi\wedge\phi_2= \partial\xi\wedge\phi & \mbox{if $\phi_1=1$}\\
(\xi\cdot\phi_1)\wedge\phi_2 & \mbox{if $\phi_1\in\Lambda^1\Phi$}\\ 0 &
\mbox{otherwise}
\end{array}
\right.
\end{equation}
and so in $Hom_k(\Lambda^*\Phi,\Xi)$, the bracket
$$[\xi,\eta]_{Hom(\Xi)}(\phi)=
(\xi\circ\overline{\hat\delta(\eta)})(\phi)- (\eta\circ\overline{\hat\delta
(\xi)})(\phi)=0$$
because the coderivations in the definition of the bracket have image in
$\Lambda^n\Phi$ with $n>0$.

It is important to note that the bracket on $Hom_k(\Lambda^*\Phi,\Xi)$ does
not restrict to the original bracket on $\Xi$ except in the abelian case;
we must introduce the "correction" term $C$.

We continue with our construction and introduce the map
$$C:\Xi\otimes\Xi\rightarrow Hom_k(\Lambda^*\Phi,\Xi)$$ by defining
$C(\xi,\eta)(\phi)=[\xi,\eta]_{\Xi}$ if $\phi=1$ and $0$ otherwise. Here,
$[\cdot,\cdot]_{\Xi}$ is the original Lie bracket on $\Xi$. Next, we must
check that $\hat\delta
C(\xi,\eta)=[\hat\delta(\xi),\hat\delta(\eta)]_{Hom(\Phi)}$ (notation as
follows equation (1)).

So for $\phi\in\Lambda^*\Phi$, we have

\begin{equation}
\hat\delta C(\xi,\eta)(\phi)=
\begin{cases}
\partial [\xi,\eta]_{\Xi} & \text{if $\phi =1$}\\ \	[ \xi,\eta]_{\Xi}
\cdot\phi & \text{if
$\phi\in\Lambda^1\Phi$} \\
0 & \text{otherwise.}
\end{cases}
\end{equation}

On the other hand, we have
$$ [\hat\delta(\xi),\hat\delta(\eta)]_{Hom(\Phi)}(\phi)=
(\hat\delta(\xi)\circ\overline{\hat\delta(\eta)}-\hat\delta(\eta)\circ
\overline{\hat\delta(\xi)})(\phi)$$
\begin{equation} =
\begin{cases}
\hat\delta(\xi)(\partial\eta\wedge\phi_2)-\hat\delta(\eta)(\partial\xi
\wedge\phi_2) & \text{if $\phi_1=1$}\\
\hat\delta(\xi)((\eta\cdot\phi_1)\wedge\phi_2)-\hat\delta(\eta)((\xi\cdot
\phi_1)\wedge\phi_2) &
\text{if $\phi_1\in\Lambda^1\Phi$}\\ 0 & \text{otherwise.}
\end{cases}
\end{equation}

The first term is non-zero only if $\phi_2=1$ in which case
$\phi=1$ and
we have $\xi\cdot(\partial\eta)-\eta\cdot(\partial\xi)$ which is equal to
$\partial[\xi,\eta]_{\Xi}$ by our original assumption on $\partial$. The
second
term is non-zero only for $\phi_2=1$ and $\phi_1\in \Lambda^1\Phi$ and is
then equal to $\xi\cdot(\eta\cdot\phi)-\eta\cdot(\xi\cdot\phi)$ which in
turn is equal to $[\xi,\eta]\cdot\phi$ by the Lie module action of $\Xi$ on
$\Phi$. Thus the BBvD hypothesis is satisfied.

Now we apply our Theorem 2 above to impose an sh-Lie structure on the
graded vector space $\VV=\{\VV_n\}$ with $\VV_0=\Xi,\VV_1=\Phi$
and $\VV_n=0$ otherwise.

It is easy to see that our construction
gives back the usual Lie algebra structure on the graded vector space $\VV,$
the semi-direct product of the Lie algebra $\Xi$ and the module $\Phi.$

\section{On shell gauge symmetries}

Up to this point we have focused primarily on unravelling the algebraic
structure implicit in the BBvD
hypothesis. This hypothesis is {\it trivially}
satisfied for classical physical theories such as general relativity and
Yang-Mills theories in the sense
that the gauge symmetries of these
physical theories satisfy the strict Lie version discussed in section 6.
On the other hand, the BBvD hypothesis appears to be precisely the condition
satisfied by the symmetries of ``free differential algebras'' which are useful
in a careful description of the Sohnius-West model of supergravity, see for
example \cite{sgrav}
and \cite{fda}. (Physicists refer to ``free differential algebras'' meaning
differential graded commutative
algebras which are free as graded commutative algebras.) Note that the
latter paper shows
that ``free differential algebras'' satisfy the BBvD hypothesis
(see equation 4.16 in \cite{fda}) without any
extra terms that vanish on shell.
Consequently,
some analysis such as the one developed in section 5 is required for a full
understanding of the algebraic
structure of these transformations.

Field dependent gauge symmetries appear in other field theories as well,
including
the class due to Ikeda \cite{ikeda} and Schaller and Strobl
\cite{schaller-strobl} and
employed by Cattaneo and Felder \cite{cf:defquant} to implement
Kontsevich's deformation quantization
\cite{kont:defquant} referred to above. These field symmetries do not
satisfy the BBvD hypothesis as we have
described it above, but rather satisfy the BBvD hypothesis``on shell". In
this section we outline how our work
may be generalized so that in the next section we can show how to apply it
to such field theories,
illustrating this in terms of one due to Ikeda (and also that of Cattaneo
and Felder).

First we explain what is meant when one says that a condition holds ``on
shell".
In essence one means that the condition holds not for all the fields of the
physical
theory, but rather that it holds only for those fields which satisfy the field
equations. In all the theories of interest here, the field equations are
Euler-Lagrange equations.
Such equations are obtained from the Lagrangian of the
physical theory. In our case we assume that the Lagrangian is a polynomial
in the
components of both the fields and their derivatives. These components may
be regarded as smooth functions on the space-time manifold $M$ and
consequently the Lagrangian is a mapping from the space $\Phi_0$ of
physical fields into $C^{\infty}M$ such that $$L(\phi)={\cal
P}_L(\phi^a,\partial_I\phi^a)$$ where ${\cal P}_L(u^a,u^a_I)$ is a
polynomial over $C^{\infty}M$ in the indeterminants $u^a,u^a_I $ and where
$\phi^a$ are the components of a typical field
in $\Phi_0$ ($I$ is a symmetric multi-index). In that which follows,
we identify $u^a$ with
$u^a_I$ where $I$ is empty. Similarly $\phi^a=\partial_I \phi^a$ where $I$
is empty. The
``action" of the physical theory is then the integral of the Lagrangian
over the space-time manifold
$M.$ All of the theories discussed in \BBvDc the supergravity example
mentioned above, and the example
due to Ikeda, discussed more fully below, are polynomial Lagrangian field
theories in the
sense we have
described above.
The Euler operator $E_a$ applied to the Lagrangian $L$ produces the
Euler-Lagrange differential operator
$E_aL$ which acts on fields via
$$E_aL(\phi)= (-1)^{|I|} \quad \partial_I (\frac{\partial {\cal P}_L}{\partial
u^a_I}(\partial_I \phi^b)).$$
Since the Lagrangian is polynomial in the components $\{\phi^a\}$ of the
fields and their derivatives
$\{\partial_I \phi^a\},$ the Euler-Lagrange differential operator is also a
mapping from $\Phi_0$ into
$C^{\infty}M$
which factors
through an appropriate polynomial over $C^{\infty}M.$

Observe that
each homogeneous polynomial ${\cal P}(u^a_I)$ of degree $k$ uniquely
defines a symmetric multi-linear mapping $\beta$ from ${\cal U}_1\times
{\cal U}_2 \times \cdots \times {\cal U}_k$ into polynomials in
$\bigcup_{i=1}^k{\cal U}_i$ such that $${\cal
P}(u^a_I)=\beta(u^a_I,u^a_I,\cdots,u^a_I)$$ for appropriate indeterminates
${\cal U}_i=\{u(i)^{a_i}_{I}\}.$ The polynomial ${\cal P}_L$ is a sum of
homogeneous terms,
each of which can be recovered from an appropriate symmetric multi-linear
mapping by evaluating the multi-linear mapping on the diagonal.

Consequently, each Lagrangian
$L$ uniquely identifies an element $\beta^L=\sum_i \beta^L_i,$ where $
\beta^L_i\in
Hom(\wedge^i_{C^{\infty}M}\partial \Phi_0,C^{\infty}M),$ such that
$$L(\phi)=\sum_i\beta^L_i(\partial_I \phi^a,\partial_I\phi^a,\cdots,
\partial_I\phi^a),$$
where $\partial \Phi_0$ denotes the vector space of the components of the fields 
and their derivatives.
We refer to this identification as {\em polarization} and will be more
precise in our algebraic formulation below. Similarly the Euler-Lagrange
differential operator admits an analogous polarization.

It is probably useful to establish a dictionary relating our algebraic approach
to field theory to more usual approaches. The algebra ${\cal A}$ is identified
with the algebra $C^{\infty}M$ of smooth functions on the space-time M and
$\Phi_0$
with the space of all physical fields of the theory. This space of fields in
simple cases is the space of all maps from $M$ into a finite-dimensional vector
space $W.$ The module $\Phi$ is an algebraic way of formulating ``jets" of
fields and $\partial$ is the map which assigns the jet $\partial
\phi=(\partial_I\phi^a)e_a^I$ to the field $\phi=\phi^ae_a.$
Elements of $Hom(\wedge ^*_{\cal A}\Phi,{\cal A})$ are identified with
``polynomials" in the fields.

In our algebraic formulation, we let
${\cal A}$ denote any commutative associative algebra and let $\Phi_0$
denote an
arbitrary
${\cal A}$-module freely and finitely generated over ${\cal A}$ with basis
$\{e_a\}.$ Working locally, we assume the existence of a finite number of
derivations $\partial_{\mu}$ of ${\cal A}$ which
admit extensions as ${\cal A}$-derivations of $\Phi_0$ in the sense that
for each ${\mu},$
$\partial_{\mu}(e_a)=0$ and
$\partial_{\mu}(f\phi)=f\partial_{\mu}\phi+(\partial_{\mu}f)\phi,$ for
$f\in {\cal A},\phi\in \Phi_0.$
For each symmetric multi-index $I=(i_1,i_2,\cdots,i_k),$ let
$\partial_I=\partial_{i_1}\circ \cdots \circ \partial_{i_k}$ and let $\Phi$
denote the ${\cal A}$-module freely generated by symbols $\{e^a_I\}$ so that
$$ \Phi=\{\phi^a_Ie_a^I | \phi_a^I\in {\cal A} \}.$$ In this context,
$L$ and the Euler-Lagrange differential operators $E_aL$ are identified
with their polarizations which are special elements of $Hom(\wedge ^*_{\cal
A} \Phi,{\cal A})$ where $\wedge^*_{\cal A} \Phi$ is the free nilpotent
cocommutative coalgebra generated by $\Phi$ over $\cal A.$

More precisely, when
we say that $L:\Phi_0 \longrightarrow {\cal A}$ is a polynomial Lagrangian,
we mean that there is a unique
$\beta^L\in Hom(\wedge_{\cal A}\Phi,{\cal A})$ such that
$$L(\phi)=\sum_i\beta^L_i(\partial \phi,\partial \phi,\cdots,\partial \phi),$$
where $\beta^L_i\in Hom(\wedge ^i_{\cal A} \Phi,{\cal A})$ is homogeneous and
$\partial$ is the mapping from $\Phi_0$ into $\Phi$ defined by
$\partial \phi= \partial_I\phi^ae_a^I.$

Here, of course, we mean that $\beta^L=\sum_i\beta^L_i$ where, for each
$i,$ $\beta^L_i$ can
only be nonzero on $\wedge^i_{\cal A}\Phi,$
i.e., for each $i, \beta^L_i$ is multilinear and symmetric having the property
that when it is evaluated on $\partial \phi\wedge \cdots \wedge \partial
\phi$ one obtains
precisely that term in ${\cal P}_L(\phi^a,\partial_I\phi^a)$ of degree $i.$
To obtain all
the terms in
$L(\phi),$ one must sum over all the homogeneous terms which appear in the
polynomial
${\cal P}_L$ which determines $L.$
It is possible to recover the mapping
$E_aL:\Phi_0 \longrightarrow {\cal A}$ from an element of $Hom(\wedge_{\cal
A}\Phi,{\cal A})$
in a similar manner. Consequently, in that which follows, we identify $L$
with $\beta^L$ and
we regard both $L$ and $E_aL$ as elements of $Hom(\wedge_{\cal A}\Phi,{\cal
A}).$

In this formulation, the ``shell"
is the subset
$\Sigma$ of $\Phi_0$ defined by
$$\Sigma=\{ \phi\in \Phi_0 |\quad E_aL(diag(\partial \phi))=0 \},$$ where
the diagonal mapping
$diag:\Phi\longrightarrow \wedge^*\Phi$ is defined by $$ diag(\phi)=\sum_p\
{\phi}^p,$$ and where
$$ {\phi}^p=
(\phi\wedge \phi\wedge\cdots\wedge \phi)\in\wedge ^p_{\cal A}
\Phi.$$ It is required that
$E_aL\in Hom(\wedge^*\Phi,{\cal A})$ be zero only on the diagonal as the
restriction of
$E_aL$ to the diagonal agrees with the polynomial counterpart of $E_aL$  and
it is the zero set of this latter function
which defines the solution space of the usual Euler-Lagrange operator. Note
that
$\Sigma$ is {\it not} a subspace of $\wedge^*\Phi.$

Define a subspace ${\cal I}$ of $Hom(\wedge^*\Phi,{\cal A})$ by $${\cal
I}=\{f\in Hom(\wedge^*\Phi,{\cal A}) | f(diag(\partial \phi))=0,
\phi\in \Sigma \}.$$
Similarly, define a subspace ${\cal N}$ of $Hom(\wedge^*\Phi,\Phi)$ by
$${\cal N}=\{ \nu\in Hom(\wedge^*\Phi,\Phi)| \nu(diag(\partial \phi))=0,
\phi \in \Sigma \}.$$
We say that $f\in Hom(\wedge^*\Phi,{\cal A})$ and $\nu\in
Hom(\wedge^*\Phi,\Phi)$ vanish ``on shell" iff $f$ and $\nu$ are in ${\cal
I}$
and ${\cal N},$ respectively.

Elements of
${\cal I}$ are ``polynomials" such as $E_aL$ which vanish ``on shell". The
``polynomials" referred to here are actually mappings from $\Phi_0$ to
${\cal A}$
which factor through polynomials over ${\cal A}$ in the indeterminates
$\{u^a_I\}$ as in
our description of the Lagrangian
$L$ above.
$Hom(\wedge ^*_{\cal A}\Phi,\Phi)$ plays the role of vector fields with
coefficients from
$Hom(\wedge^*_{\cal A}\Phi,{\cal A}),$ and ${\cal N}$ plays the role of the
space of
vector fields whose coefficients vanish on $\Sigma.$

At this point, we generalize the BBvD hypothesis as follows. We say that
the k-linear
mapping $\delta: \Xi\longrightarrow Hom(\wedge ^*_{\cal A}\Phi,\Phi)$
satisfies the
{\it generalized \BBvD hypothesis}, denoted gBBvD, iff there exists a
skew-symmetric k-bilinear mapping $C:\Xi \times \Xi \longrightarrow
Hom(\wedge ^*_{\cal A}\Phi,\Xi)$ and an extension
$\hat\delta$ of $\delta$ to $Hom(\wedge ^*_{\cal A}\Phi,\Xi)$ such that
$$[\delta(\xi),\delta(\eta)]-\hat \delta (C(\xi,\eta))\in {\cal N}.$$ for
all $\xi,\eta\in \Xi.$ Thus the BBvD hypothesis of Section 3 holds ``on
shell".
A consequence
of this hypothesis is that there exists a skew-symmetric mapping $\nu:\Xi
\times \Xi
\longrightarrow {\cal N}$ such that
$$[\delta(\xi),\delta(\eta)]=\hat \delta (C(\xi,\eta))+\nu(\xi,\eta).$$

Utilizing this mapping $C,$ one can define a bracket on $\LPX$ analogous to
that defined before in the presence of the BBvD hypothesis:
$$
[\pi_1,\pi_2] := \pi_1\odot \hat\dd(\pi_2) - \pi_2\odot \hat\dd(\pi_1) +
C(\pi_1,\pi_2).
$$

Injectivity of $\hat\delta$ is not easily obtained and seems to be needed
to obtain
a proof of the Jacobi identity. Thus, in general the bracket on $\LPX$ will not
satisfy
the Jacobi identity.

On the other hand, we can use the calculations of the proof of
Theorem 1 to show that
$$[\hat\delta(\pi_1),\hat\delta(\pi_2)]-\hat\delta([\pi_1,\pi_2])\in {\cal N}$$
for all $\pi_1,\pi_2$ in $\LPX.$
By using the calculations in the proof of Theorem 2, it is easy to show that:
$$D\circ \overline D(\xi_1\wedge \xi_2\wedge\phi_1\wedge \cdots \wedge
\phi_n)=$$
$$= \ -[\delta(\xi_1),\delta(\xi_2)](\phi_1\wedge \cdots \wedge \phi_n)+
\hat\delta(C(\xi_1,\xi_2))(\phi_1\wedge \cdots \wedge \phi_n)$$ $$ =\
\nu(\xi_1,\xi_2)(\phi_1\wedge \cdots \wedge\phi_n)$$ and that
$$D\circ \overline D(\xi_1\wedge \xi_2 \wedge \xi_3\wedge\phi_1 \wedge
\cdots \wedge \phi_n)=
Jacobi(\xi_1,\xi_2,\xi_3)(\phi_1\wedge \cdots \wedge \phi_n)$$ where
$$Jacobi(\xi_1,\xi_2,\xi_3)=
([[\xi_1,\xi_2],\xi_3]-[[\xi_1,\xi_3],\xi_2]+[[\xi_2,\xi_3],\xi_1]).$$ In
this latter equation we have used the notation $[\xi,\eta]$ in place of
$C(\xi,\eta).$

{\it It follows from these equations that $D\circ \overline D$ is zero ``on
shell" provided that both
the generalized BBvD hypothesis holds and that $Jacobi(\xi,\eta,\zeta)$ is
zero on shell for
arbitrary constants $\xi,\eta,\zeta\in Hom(\wedge ^*_{\cal A}\Phi,\Xi).$}

In the example due to Ikeda \cite{ikeda}
discussed in detail in Section 8, it is easy to prove that
$Jacobi(\xi_1,\xi_2,\xi_3)=0$ (not just zero ``on shell") using the
equation immediately
prior to equation 2.10 in his paper.
Consequently the gauge symmetries of this Poisson $\sigma$-model satisfy
the postulates
of an sh-Lie algebra ``on shell".

\section{A $\Sigma$-model example}
In Ikeda's paper \cite{ikeda}, there is a finite dimensional vector space $V$
with basis $\{T_A\}$ which later we will show is the dual of a Poisson
manifold.
We do this via a generalization of the classical Kirillov-Kostant bracket
which exhibits
the dual $\mathfrak g^*$ of a Lie algebra $\mathfrak g$ as a Poisson manifold.

In our analysis of Ikeda's example, our space $\Xi$ is the space of maps
$Maps(\Sigma,V)$
and the space $\Phi_0$
is the set of ordered pairs $\phi=(\psi,h)$ where:

(1) $\psi$ is a mapping from a given two-dimensional manifold $\Sigma$ into the
dual $V^*$ of the vector space $V,$ and

(2) $h$ is a mapping from the same manifold $\Sigma$ to $T^*\Sigma\otimes
V,$ which in fact is required to be a section of the vector bundle
$T^*\Sigma\otimes V \longrightarrow \Sigma.$ These mappings are denoted
locally by
$\psi(x)=\psi_A(x)T^A$ and $h(x)=h_{\mu}^A(dx^{\mu}\otimes T_A),$ where
$\{T_A\}$ is
a basis of $V$ and $\{T^A\}$ is the basis of $V^*$ dual to $\{T_A\}.$

For the most part, our exposition follows that of Ikeda, although we use
the notation
$\phi= (\psi,h)$ for the fields of the theory whereas Ikeda's notation for
the fields is $(\phi,h).$ We also denote Ikeda's vector space $M$ by $V.$
As is the case earlier in the paper, the space $\Phi$ denotes the ${\cal
A}= C^{\infty}M$
module whose elements are $\phi^a_Ie_a^I$ where $\{e^I_a\}$ is a basis of
the module and $\phi^a_I\in
{\cal A}.$ This formulation is our algebraic description of the jet bundle
of the vector bundle whose
sections are the fields
$\Phi_0.$ Ikeda would denote $\phi^a_I$ as $\partial_I\phi^a.$ Observe that
$\Phi_0$ may be identified as a subspace of $\Phi.$

There
is a parallel development to Ikeda's work in Cattaneo and Felder
\cite{cf:defquant} in
which $\Sigma$
is a 2-dimensional disc and the target
(denoted by $M$ in Cattaneo and Felder) is an arbitrary Poisson manifold.
It is not
hard to see that the ordered
pairs $(\psi,h)$ of Ikeda may in fact be interpreted in a manner similar to
that in the exposition of
Cattaneo and Felder where $\psi:\Sigma \longrightarrow M$ is an arbitrary
smooth mapping ($\psi$ is denoted by $X$ in Cattaneo and Felder) and $h$ is
a section
(denoted by $\eta$ in Cattaneo and Felder) of the bundle
$\psi^*(T^*M)\otimes T^*\Sigma \longrightarrow \Sigma$ (notice that the
factors in their tensor product are reversed from the conventions used in
our description of Ikedas' results). In their exposition the section $h$
may be
written as $h(x)=h_{i,\mu}(x)(dx^i\otimes du^{\mu})$ where $\{dx^i\}$ is a
basis of $T^*_{\psi(x)}M,$
which, in the case $M$ is a vector space $V,$ may be identified with a fixed
basis $\{T^A\}$ of $T^*_0V=V^*.$

When one compares these two approaches, one sees that Ikeda's target space
is the vector space we
have called $V^*$ while $\Sigma$ is an arbitrary 2-dimensional manifold,
whereas for Cattaneo and
Felder
$\Sigma$ is a disc $D$ and the target space $M$ is a general Poisson
manifold. The parallel between the two is closer than one might initially
expect since Ikeda uses the vector space $V$ to generate a
Poisson structure on $V^*.$

Ikeda proceeds to investigate possible
gauge symmetries $\delta(c)$ {\em before} looking for Lagrangians. The
gauge symmetry mapping
$\delta$ is defined locally, in this theory, as follows.
Let ${\mathcal P}$ denote the commutative polynomial algebra generated by
the basis $\{T_A\}.$
Let $\pi_A,
\pi^A_{\mu}, \pi^A$ denote the projections defined by
$\pi_A(\phi)=\pi_A(\psi,h)=\pi_A(\psi)=\psi_A$,
$\pi^A_{\mu}(\phi)=\pi^A_{\mu}(\psi,h)=\pi^A_{\mu}(h)=h^A_{\mu}$ and
$\pi^A(c)= c^A,$ respectively. Consider arbitrary polynomials $\{W_{AB}\}$
in ${\cal P}$ and define the components of $\delta(c)(\phi)$ by
$$\pi^A_{\mu}(\delta(c)(\phi))=\partial_{\mu}c^A+ \frac {\partial
W_{BD}(\psi)} {\partial T_A} h_{\mu}^Bc^D$$ and
$$\pi_A(\delta(c)(\phi))=W_{BA}(\psi)c^B.$$

Here $W_{AB}(\psi)$ is a concise notation for the polynomial $W_{AB}$
evaluated by replacing the
generators $\{T_A\}$ by the correponding components $\{\psi_A\}$ of $\psi,
$ that is, $W_{AB} = W_{AB}^a T_a$ and $W_{AB}(\psi) = W_{AB}^a \psi_a$
where $a$ is a symmetric multi-index,
$\psi_a=\psi_{A_1}\psi_{A_2}\cdots\psi_{A_n}$ and similarly for $T_a.$

Notice that, in case $V$ is a Lie algebra and $W_{AB}(T)=[T_A,T_B]=
f^C_{AB}T_C,$ the polynomials
$W_{AB}$ define the Lie algebra structure on the vector space $V$ with
structure
constants $\{f^C_{AB}\}.$
This then induces a Lie algebra structure on the parameter space $\Xi$ of
all mappings $c$ from
$\Sigma$ into
$V,$ as one expects in traditional Yang-Mills theory. In this case, the
$\psi$-component of $\delta(c)$ is the
coadjoint action of the parameter space $\Xi$ on the space of maps from
$\Sigma$
into $V,$ while the
$h$-component is simply the ``covariant derivative" of $c$ relative to the
connection defined by the gauge field
$h.$
Thus, by introducing more general polynomials $W_{AB},$ Ikeda is
introducing a generalization
of ordinary gauge theory by requiring that the gauge symmetries be defined
via the polynomials
$W_{AB}.$ For this
generalization to work, Ikeda imposes restrictions on the polynomials
$W_{AB}$ which amount to making
$\cal P$
a Lie algebra, hence his terminology of `non-linear Lie algebra'. In order
to obtain an algebraic
structure on ${\cal P}$ analogous to the usual Lie structure required in
gauge theory, Ikeda's bracket is defined on generators of $\cal
P$ by
$$[T_A,T_B] = W_{AB}\in \cal P$$
and extended to all of $\cal P$ via the Leibniz rule: $[T_A,\ ]$ and $[\
,T_B]$ are derivations
of the commutative algebra $\cal P$. Ikeda requires that these polynomials
satisfy conditions which make
$\cal P$ a Poisson algebra. Thus the polynomials $\{W_{AB}\}$ in ${\cal P}$
are subject to skew-symmetry: $W_{AB}=-W_{BA}$ and an appropriate
generalization of the
usual coordinate form of the Jacobi condition:

$$W_{AD}\frac {\partial W_{BC}}{\partial T_D}+W_{BD}\frac {\partial
W_{CA}}{\partial T_D}+
W_{CD}\frac {\partial W_{AB}}{\partial T_D}=0 .$$

To see $V^*$ as a Poisson manifold, we will imbed $V$ in $V^{**}$ as the
linear functionals and thus regard the algebra ${\cal P}$ as the subalgebra
of $C^\infty (V^*)$ generated by the basis $\{T_A\}.$

Regarding $T_A$'s as functions on $V^*$, we have a bi-vector field
$$W_{AB}\frac {\partial }{\partial T_A}\wedge \frac {\partial }{\partial
T_B}$$ on $V^*.$

This makes $V^*$ a Poisson
manifold with
$$\{f,g\} :=
W_{AB}\frac {\partial f}{\partial T_A}\wedge \frac {\partial g}{\partial
T_B}$$ \noindent
for $f,g\in C^\infty (V^*)$.
Now notice that, using Ikeda's notation as defined above, we have for each
$c$ that $\delta(c)$ is a
mapping from
$\Phi_0$ to $\Phi_0.$ Since $\Phi_0$ is a vector space, it follows that
with any reasonable topology
on $\Phi_0$ one can identify the tangent space of $\Phi_0$ at a point $\phi\in
\Phi_0$ with $\Phi_0$ itself. Thus maps from $\Phi_0$ into $\Phi_0$ may be
regarded as vector fields on $\Phi_0.$
Recall that $\delta(c)$ is a vector field on the space $\Phi_0$
of fields. Thus $\delta(c)(\phi)$ is a tangent vector to $\Phi_0$
at $\phi.$

By an obvious abuse of notation one may write:
\begin{equation} \label{vf} \delta(c)(\phi)=(W_{BA}(\psi)c^B)\frac
{\partial}{\partial
\psi_A}+ (\partial_{\mu}c^A+
\frac {\partial W_{BD}(\psi)} {\partial T_A} h_{\mu}^Bc^D) \frac
{\partial}{\partial
h^A_{\mu}}.\end{equation}
\noindent

Now the components of $\delta(c)(\phi)$ as defined in equation (8)
are polynomials in the
components of $\phi.$ Consequently, in conformity with our conventions in
section 7, we can identify
$\delta(c)$ with the unique element of
$Hom(\wedge^*\Phi,\Phi)$ whose value at $diag(\partial \phi),$ for $\phi\in
\Phi_0,$ gives
$\delta(c)(\phi)$ as defined by Ikeda.

The usual Lie bracket of the
vector fields $\delta(c_1)$ and $\delta(c_2)$ as defined by equation
\ref{vf} corresponds to our Lie
structure on $Hom(\wedge^*\Phi,\Phi).$
Using his brackets, Ikeda finds that the $\psi$ component of
$[\delta(c_1),\delta(c_2)](\phi)$
is given by $$[\delta(c_1),\delta(c_2)](\psi)=\delta(c_3(\psi))(\psi)$$ where

\begin{equation} \label{comm}\pi_A(c_3(\psi)) =\frac {\partial
W_{BD}}{\partial T_A}(\psi)c_1^Bc_2^D.
\end{equation}

We see that the Lie bracket of $[\delta(c_1),\delta(c_2)]$ is not of the
form $\delta(c)$
where $c_3$ is a gauge parameter independent of the fields $\phi$ but
rather the gauge
parameter $c$ depends on $c_1,c_2$ and on the field $\psi.$ Thus one does
not have closure
on the original space of gauge parameters $Maps(\Sigma, V)=\Xi.$ We are forced
to enlarge the space of gauge parameters $\Xi$ to include mappings from
$\Phi$ to $\Xi$.

In Ikeda's context, these mappings are polynomials in the components of the
fields and their
derivatives. Consequently, they are identified with elements of $\LPX$ in
our formulation.

If we had only the fields $\psi$ to deal with, the BBvD hypothesis would be
satisfied and
we would be able to apply the ideas in earlier sections to describe Ikeda's
algebra of gauge
transformations
as an sh-Lie algebra on a graded space with $\Xi$ in degree zero. However,
the $h$-component
transforms more
subtly. To handle this, Ikeda makes the definition \begin{equation}
\label{covd}
D_{\mu}\psi_A=\partial_{\mu}\psi_A+W_{AB}(\psi)h^B_{\mu}. \end{equation}
\noindent
(The resemblance to a
covariant derivative is formal; it is not yet understood as arising in an
obvious manner from a ``representation"
of the nonlinear Lie algebra defined by Ikeda.) He then calculates

\begin{equation} [\delta(c_1),\delta(c_2)](h)=\delta(c_3(\psi))(h)- (\frac
{\partial^2W_{CD}}{\partial \psi_A \partial
\psi_B}(D_{\mu}\psi_B)c_1^Cc_2^D) \frac
{\partial}{\partial h^A_{\mu}}.
\label{offshell}
\end{equation}
\noindent
Thus the BBvD hypothesis fails, but
the generalized hypothesis holds.
When closure on the
original space of $V$-valued gauge parameters is lost, physicists speak of
an `open algebra'.

Having established his gauge algebra and potential gauge symmetries, Ikeda
then searches for an
appropriate Lagrangian. Up to a total divergence, the Lagrangian of Ikeda's
theory is
$${\cal L}=\epsilon^{\mu\nu}\{h^A_{\mu}D_{\nu}\psi_A- \frac{1}{2}
W_{AB}(\psi)h^A_{\mu}h^B_{\nu}\}.$$
\noindent This includes self-interacting terms for the generalized gauge
fields $h$ along with a
minimal coupling of the scalar field $\psi$ through the generalized
covariant derivative defined in
equation (10)  above. The tensor $\epsilon^{\mu\nu}$ is the area
element which is assumed to be present on $\Sigma.$

Ikeda really
works with an equivalent Lagrangian which differs from the one given above
by a divergence,
although the physical content of the Lagrangian defined above is clearer.

Ikeda shows that for his equivalent
Lagrangian $${\cal L}(\phi,\partial \phi)= {\cal L}(\psi,h,\partial
\psi,\partial h),$$ the function
$\delta(c)({\cal L})$
is a divergence for all parameters
$c.$ This is precisely the property physicists require in order to call
$\delta$ a gauge
symmetry.

The field equations of the Lagrangian are $$D_{\mu}\psi_A=0 \quad \quad
R^A_{\mu\nu}=0$$
where $R^A_{\mu\nu}$ is the ``generalized" curvature
$$R^A_{\mu\nu}=\partial _{\mu}h^A_{\nu}-\partial_{\nu}h^A_{\mu}+
\frac{\partial W_{BC}}{\partial T_A}(\psi)h^B_{\mu}h^C_{\nu}$$ of the
``generalized gauge field" $h=h^A_{\mu}(dx^{\mu}\otimes T_A).$

The coefficient of the last term on the right-hand side of equation
(\ref{offshell}) is polynomial
in the components of $\phi\in \Phi_0$ and their derivatives and (as in
section 7) determines a unique bilinear mapping $\nu$ from
$\Xi\times\Xi$ into ${\cal N}=\{ \nu\in Hom(\wedge^*\Phi,\Phi)|
\nu(diag(\partial \phi))=0, \phi \in \Sigma \}$ such that \begin{equation}
\label{modI}
[\delta(c_1),\delta(c_2)]=\delta(C(c_1,c_2))+\nu(c_1,c_2) \end{equation}
\noindent
where
$C(c_1,c_2)=c_3:\wedge^*\Phi\longrightarrow \Xi$ is defined by equation
(\ref{comm}).
This latter property (12)
is the one we have referred to above as the gBBvD hypothesis.

A similar analysis applies to the Lagrangian of Cattaneo and Felder. Thus
Ikeda and Cattaneo
and Felder provide examples of field theories which satisfy the generalized
BBvD hypothesis and it
is this condition which we have assumed in sections 7 and 8. The gauge
symmetries of these theories require
a modification of the sh-Lie structure one obtains from the gauge
structures of field theories satisfying the BBvD hypothesis.

\providecommand{\bysame}{\leavevmode\hbox to3em{\hrulefill}\thinspace}

\end{document}